\documentclass[12pt, reqno]{amsart}

\usepackage{amssymb,latexsym,amsmath,amsfonts}
\usepackage{mathrsfs}
\usepackage{graphicx}
\usepackage[usenames]{color}
\usepackage{hyperref}
\usepackage{comment}
\usepackage{enumitem}

\definecolor{DPurple}{rgb}{0.46,0.2,0.69}

\numberwithin{equation}{section}

\allowdisplaybreaks

\theoremstyle{definition}
\newtheorem{definition}{Definition}[section]

\newtheorem{question}{Question}[section]
\theoremstyle{remark}
\newtheorem{remark}{Remark}[section]

 \theoremstyle{plain}
\newtheorem{theorem}[definition]{Theorem}
\newtheorem{result}{Result}[section]

\newtheorem{lemma}[definition]{Lemma}

\newtheorem{example}{Example}[section]



\begin{document}

\title{A Function--Sharing Criterion for Normal Functions}

\author{Gopal Datt}
\address{Department of Mathematics, Babasaheb Bhimrao Ambedkar University, Lucknow, India}
\email{ggopal.datt@gmail.com, gopal.du@gmail.com}
\author{Ritesh Pal}
\address{Department of Mathematics, Babasaheb Bhimrao Ambedkar University, Lucknow, India}
\email{rriteshpal@gmail.com, ritesh.rs.math@bbau.ac.in}
\author{Ashish Kumar Trivedi}
\address{ Department of Mathematics, University of Delhi, Delhi 110007}
\email{trivediashish2016@gmail.com, aktrivedi@maths.du.ac.in}


\maketitle
\begin{abstract}

In this paper, we present a function\,-\,sharing criterion for the
normality of meromorphic functions. Let $f$ be a meromorphic function
in the unit disc $\mathbb{D}\subset \mathbb{C}$, $\psi_1$, $\psi_2$, and
$\psi_3$ be three meromorphic functions in the unit disc \(\mathbb{D} \),
continuous on $ \partial{\mathbb{D}}:=\{z\in\mathbb{C}\,:\,|z|=1\}$, such that
$\psi_i(z)\neq\psi_j(z)$ $(1\leq i<j\leq 3)$ on $\partial\mathbb{D}$.
We prove that, if $\psi_1$, $\psi_2$, and $\psi_3$ share the function
$f$ on $\mathbb{D}$, then $f$ is normal. Building upon this,
we further establish an additional criterion for normal functions.%
\end{abstract}

\section{Introduction and Main Results}
The study of convergence phenomena in complex analysis has a long history.
Already, in the nineteenth century, Weierstrass established that uniform
convergence preserves holomorphy, a principle that underpins much of modern function theory.
Subsequent advances by Stieltjes, Osgood, Porter, and Vitali gradually expanded
the range of situations where convergence could be guaranteed,
laying the groundwork for Montel's celebrated theory of normal families in the early twentieth century.\smallskip

Montel's insight was that boundedness and equicontinuity,
when combined with the Ascoli\,--\,Arzel\'a theorem, yield
compactness of families of holomorphic or meromorphic functions,
thereby providing a powerful method for proving theorems of Picard,
Schottky, Landau, and others. Interested readers may refer the
monograph by Schiff~\cite{Schiff1993} for vivid introduction to
normal families.\smallskip

Although Montel's terminology ``normal family'' dates from 1911, the idea
soon inspired a parallel development: to attribute the property of
normality not merely to a family, but to a single meromorphic function.
This point of view was first suggested by Yosida in 1934 in \cite{Yosida1934}
where he talked about the \emph{property\,$(A)$}: Given any sequence of complex number
$(a_{i})$ the family of meromorphic function $y_{i}=y(z+a_{i})$, where $i=1,2,...$
is normal family in any closed and finite domain. The class of meromorphic functions
having the property\,$(A)$ was called the class\,$(A)$ by Yosida.\smallskip

This was subsequently developed in detail by Noshiro in~\cite{Noshiro1938},
where he studied a class of meromorphic functions in the unit disc:
Suppose that $f$ be a meromorphic in the unit disc $\mathbb{D}$ and
consider the family $\{ f_{a}(z) \}$ formed by all the functions
\[
f_{a}(z) \equiv f\!\left( \frac{z-a}{\overline{a}z - 1} \right),
\]
where $a$ varies throughout the unit disc $\mathbb{D}$. Following the
notion of Yosida~\cite{Yosida1934}, Noshiro also said that the function
$f$ belongs to class\,$(A)$ if $\{ f_{a}(z) \}$ is a
normal family on unit disc  $\mathbb{D}$.\smallskip

In the same paper, Noshiro established the following characterization: A non-constant
meromorphic function $f$ in the unit disc $\mathbb{D}$
belongs to class\,$(A)$  if and only if there exists a finite constant $C$ such that
\begin{equation} \label{eq:ClassA}
   (1-|z|^{2})f^{\#}(z)\, \;\leq\; C, \ \mbox{where} \ f^{\#}(z) = \dfrac{|f'(z)|}{1+|f(z)|^{2}}.
\end{equation}
\smallskip

The name \emph{normal function} was given later by
Lehto and Virtanen in~\cite{LehtoVirtanen1957} in 1957,
when a large\,-\,scale study was undertaken and they showed
that the idea of normal meromorphic function is much related to
some of the most important problems of the boundary behaviour of
meromorphic functions.\smallskip

In their formulation, a meromorphic function $f$ in
a hyperbolic domain $D \subset \mathbb{C}$ is called $\emph{normal}$
if its orbit under the automorphism group of $D$,
\[
   \{\,f \circ \phi : \phi \in Aut(D)\,\},
\]
forms a normal family in $D$ with respect to the spherical metric.\smallskip

Normal meromorphic functions admit the following characterization in
terms of the spherical derivative: A non-constant meromorphic function
$f$ is \emph{normal} in a domain $D$ (which is necessarily of hyperbolic type)
if and only if there exists a finite constant $C$ such that
\begin{equation} \label{eq:normal}
   f^{\#}(z)\,|dz| \;\leq\; C\, d\alpha(z),
\end{equation}
where, $f^{\#}(z)$ is the spherical derivative of $f$,
and $d\alpha(z)$ denotes the element of length in the hyperbolic metric of $D$.\smallskip

The concept of normal functions has proven to be a natural and robust
extension of Montel's framework. It links directly to questions of
boundary behavior, Lindel\"{o}f-type principles, and value distribution, and it also connects to quasiconformal mappings; see~\cite{FletcherNicks2024,LehtoVirtanen1957,Lindelof1915,Lo2013}.\smallskip

Montel, in \cite{Montel1912}, presented his {\it fundamental normality criterion} ({\it Crit\`ere
fondametal}) for families of analytic function.
\begin{result}\cite{Montel1912}\label{T:fnt}
Let $\mathcal{F}$ be a family of meromorphic functions on a
domain $D \subseteq \mathbb{C}$, such that each function $f \in \mathcal{F}$
omit three distinct values say a,b,c in $\hat{\mathbb{C}}$,
then $\mathcal{F}$ is normal in $D$.
\end{result}

One of the central and elegant results about normal functions,
echoing the spirit of the fundamental normality criterion, is the following:

\begin{result}\label{R: omit}
Let \( f \) be a meromorphic function in the unit disc \( \mathbb{D} \).
If \( f \) omits at least three distinct values in \( \hat{\mathbb{C}} \),
then \( f \) is a normal function.
\end{result}

A natural question arises: can the three values in the fundamental
normality test, Result~\ref{T:fnt}, be replaced by three continuous functions?
This was achieved by Bargmann et al. in~\cite{BBHM}, where they established that a family of
meromorphic functions on a domain $D$ is normal if each member of the
family avoids three continuous functions in $D$, and these three functions
avoid each other in $D$.\smallskip

Recall that a meromorphic function $f$ avoids a function $g$ on a
domain $D$ if $f(z)\neq g(z)$, for all $z\in D$. A more precise definition
is as follows:\smallskip

We say that two functions \( f \) and \( g \) avoid each other uniformly if
there exists a constant \( \delta > 0 \) such that, for every point \( z \)
in their common domain, the spherical distance between
\( f(z) \) and \( g(z) \) is at least \( \delta \).\smallskip

Motivated by the work of Bargmann et al.~\cite{BBHM}, Lappan in ~\cite{Lappan2003}
subsequently extended the result to the class of normal functions
on the unit disc $\mathbb{D}$, as detailed below:
\begin{result}\cite[Theorem 3]{Lappan2003}\label{R: lappan}
Let \(\psi_1, \psi_2,\)\ and \(\psi_3 \) be three continuous functions
defined on the unit disc \( \mathbb{D} \), and suppose they avoid each other
uniformly. Furthermore, assume the family
\[
\{ \psi_j \circ \varphi : \varphi \in \Phi \}
\]
is normal in \( \mathbb{D} \), for each \( j = 1, 2, 3 \),
where \( \Phi = \{ \varphi : \mathbb{D} \to \mathbb{D} \mid \varphi \text{ is conformal} \} \).
Let \( f \) be a meromorphic function in \( \mathbb{D} \)
such that \( f(z) \neq \psi_i(z) \) for each \( i \in \{1,2,3\} \) in $\mathbb{D}$.
Then \( f \) is a normal function.
\end{result}
It was also noted by Lappan in ~\cite{Lappan2003}, the assumption that the functions
\( \psi_1, \psi_2 \), and \( \psi_3 \) avoid each other uniformly is necessary
to hold the result \ref{R: lappan}.\smallskip

However, when the functions \(\psi_1, \psi_2\), and \(\psi_3 \) are
meromorphic in the unit disc \( \mathbb{D} \) and continuous on
its boundary \(\partial{\mathbb{D}}:=\{z\in\mathbb{C}\,:\,|z|=1\} \), it suffices that they avoid
each other pointwise on \( \partial \mathbb{D} \), uniform avoidance
throughout \( \mathbb{D} \) is not necessary in this case, as
established by Xu and Qiu in~\cite{XuQiu2011}.

\begin{result}\cite[Theorem 1]{XuQiu2011}\label{avoid}
Let $f$ be a meromorphic function in the unit disc $\mathbb{D}$, $\psi_1, \psi_2$, and $\psi_3$
be three meromorphic functions in $\mathbb{D}$, continuous on $\partial{\mathbb{D}}$,
such that $\psi_i \neq \psi_j$ on $\partial\mathbb{D}$.
If $f(z)\neq  \psi_i(z)$, $(i = 1,2,3)$ in $\mathbb{D}$, then f is normal.
\end{result}

In the aforementioned results~\ref{avoid} and~\ref{R: lappan}, it is evident
that their normality criteria are based on the avoidance of a
meromorphic function. More recently, various criteria for the normality
of families of meromorphic functions have been developed,
involving the sharing of values, sets, and functions, see~\cite{Charak2016,DattKumar2015, DattKumar2016,ZengLahiri2015}.\smallskip

Motivated by these developments, {\it A natural question arises}: can
we employ sharing of function to establish normality for a meromorphic
function on the unit disc $\mathbb{D}$? Recall the definition:
Let $f$, $g$ and $h$ be three meromorphic functions in the domain $D$.
We say $f$ and $g$ \textbf{share} the function $h$ (ignoring multiplicity),
if $\{z\in D: f(z)=h(z)\}=\{z\in D: g(z)=h(z)\}$.

\begin{question} \label{Q: sharing}
What if,\textit{ we consider three meromorphic functions in the unit disc
$\mathbb{D}$ and continuous on $\partial \mathbb{D}$ such that they mutually
avoid each other on the boundary $\partial \mathbb{D}$,
and share a meromorphic function $f$ on $\mathbb{D}$.}
Can we conclude normality of the function $f$?
\end{question}

In this paper, we endeavor to address Question~\ref{Q: sharing},
and to the best of our knowledge, it is the \emph{first paper}
in the direction of function\,-\,sharing criterion for normal functions.
Next, we are going to proof the theorem which gives an
affirmative answer to the
Question~\ref{Q: sharing}.

\begin{theorem}\label{main}
Let $f$ be a meromorphic function in the unit disc $\mathbb{D}$,
$\psi_1$, $\psi_2$, and $\psi_3$ be three meromorphic functions in
$\mathbb{D}$, continuous on $\partial{\mathbb{D}}$ such that
$\psi_i(z)\neq\psi_j(z)$ $(1\leq i<j\leq 3)$ on $\partial\mathbb{D}$.
If $\psi_1$, $\psi_2$, and $\psi_3$ share the function $f$ on
$\mathbb{D}$, then $f$ is normal.
\end{theorem}

It is quite evident to see that, Result~\ref{avoid} by
Xu and Qiu may be viewed as a special case of our Theorem
~\ref{main}, corresponding to the
situation in which the function $f$ avoids each of the
functions $\psi_1$, $\psi_2$, and $\psi_3$ on $\mathbb{D}$.

\begin{example}We now provide an example to illustrate Theorem~\ref{main} .

\end{example}
Let $\psi_1(z) = 2z$, $\psi_2(z) = 3z$, $\psi_3(z) = 4z$
be three meromorphic functions on the unit disc $\mathbb{D}$, continuous
on the boundary $\partial\mathbb{D}$, and mutually avoiding
each other on $\partial\mathbb{D}$; that is,
$\psi_i(z) \neq \psi_j(z)$  \text{for all } $z \in \partial\mathbb{D}$, $i \neq j$.
Consider the function
$f(z) = z$,
which is meromorphic on $\mathbb{D}$. Observe that:
\[\{z \in \mathbb{D} : \psi_1(z) = f(z)\} = \{z \in \mathbb{D} :
\psi_2(z) = f(z)\} = \{z \in \mathbb{D} : \psi_3(z) = f(z)\} = \{0\}.\]

\noindent Thus, all three functions $\psi_1$, $\psi_2$, and $\psi_3$
share the function $f$. Here, $(1-|z|^2)f^\#(z)=\dfrac{1-|z|^2}{1+|z|^2}$ is bounded
on $\mathbb{D}$, hence $f$ is normal.

\begin{example}\label{Ex:sharp} The number
three in the statement of Theorem~\ref{main} is sharp.
\end{example}
Let $\psi_1(z)\equiv i$, $\psi_2(z)\equiv -i$, $\psi_3(z)\equiv 1$
be three meromorphic functions on the unit disc $\mathbb{D}$, each
continuous on the boundary $\partial\mathbb{D}$, and mutually avoiding
each other on $\partial\mathbb{D}$; that is,
$\psi_i(z) \neq \psi_j(z) \quad \text{for all } z \in
\partial\mathbb{D},\ i \neq j.$
Now consider the function
$f(z)=\tan(\frac{1}{1-z})$,
which is meromorphic on $\mathbb{D}$. Observe that:\smallskip

$\bullet$ $\psi_1$ and $\psi_2$ share the function $f$ because
\[\{z \in \mathbb{D} : \psi_1(z) = f(z)\} = \{z \in \mathbb{D}
: \psi_2(z) = f(z)\} = \varnothing.\]

$\bullet$
$\{z \in \mathbb{D} : \psi_3(z) = f(z)\} = \left\{ 1-\frac{4}{(4k+1)\pi}: k \in \mathbb{N}\cup\{0\}\right\}$.\smallskip

\noindent The three functions $\psi_1$, $\psi_2$, and $\psi_3$ do not
simultaneously share the function $f$, nor does $f$ avoid all three.\smallskip

\noindent However, the function $f(z)=\tan(\frac{1}{1-z})$ is not normal. To see this we consider
\[
(1-|z|^2)f^\#(z) = (1-|z|^2)\frac{|\sec^2(1/(1-z))|}{1 + |\tan(1/(1-z))|^2}\left(\frac{1}{|1-z|^2}\right).
\]
Which approaches towards $\infty$, when $z$ approaches towards $1$, through real axis i.e.,
\(z=r \to 1^- \Rightarrow (1-|r|^2)f^\#(r) \to \infty \).\smallskip

Now one can ask the reasonable question, can we relax the hypothesis in
Theorem~\ref{main}, so that only two functions say,
$\psi_1$ and $\psi_2$ share the function $f$ on $\mathbb{D}$?
If so, then what condition should be imposed on $\psi_3$ in order to ensure that
the function $f$ is normal on $\mathbb{D}$?
Next theorem, which we present gives one of the answer
to the question above by imposing the condition that $f$
should avoid the function $\psi_3$ on $\mathbb{D}$.

\begin{theorem}\label{submain}
Let $f$ be a meromorphic function in the unit disc $\mathbb{D}$,
$\psi_1$, $\psi_2$, and $\psi_3$ be three meromorphic functions in
$\mathbb{D}$ and continuous on $\partial{\mathbb{D}}$ such that
$\psi_i(z)\neq\psi_j(z)$ $(1\leq i<j\leq 3)$ on $\partial\mathbb{D}$.
If $\psi_1$ and $\psi_2$ share the function $f$ on $\mathbb{D}$,
and $\psi_3$ avoids the function $f$ on $\mathbb{D}$, then $f$ is normal.
\end{theorem}

In fact, Theorem~\ref{submain} also encompasses the Result~\ref{avoid}
as a special case, when it is assumed that the
function $f$ avoids each of the functions $\psi_1$ and $\psi_2$
on the unit disc $\mathbb{D}$.

\begin{example}To illustrate Theorem~\ref{submain}, we provide the following example.

\end{example}

Let $\psi_1(z) = 2z$, $\psi_2(z) = 3z$, $\psi_3(z) = 1/z$ be three
meromorphic functions defined on the unit disc $\mathbb{D}$,
continuous on the boundary $\partial\mathbb{D}$, and mutually
avoiding each other on $\partial\mathbb{D}$; that is,
$\psi_i(z) \neq \psi_j(z)$  \text{for all } $z \in \partial\mathbb{D}$,
$i \neq j$. Consider the function $f(z) = z$, which is meromorphic on
$\mathbb{D}$. We observe the following:\smallskip

$\bullet$  The functions $\psi_1$ and $\psi_2$ share the function $f$,
\[\{z \in \mathbb{D} : \psi_1(z) = f(z)\} = \{z \in \mathbb{D} :
\psi_2(z) = f(z)\} = \{0\}.\]

$\bullet$ The function $\psi_3(z) = 1/z$ avoids $f(z)=z$ on $\mathbb{D}$, because
\[\{z \in \mathbb{D} : \psi_3(z) = f(z)\} = \varnothing.\]

\noindent Here, the function $f(z)=z$ is normal.

\begin{remark}
The Example~\ref{Ex:sharp} also shows, we cannot omit the condition
that $f$ avoids $\psi_3$ in Theorem~\ref{submain}.
\end{remark}

\section{Essential lemmas}

In order to establish our main theorems, we begin by recalling
several classical lemmas that will serve as essential tools in our proofs.
The first among these is the well-celebrated Picard's Little Theorem.

\begin{lemma}\cite{Picard-1}
A non-constant entire function
takes every complex value with at most one
possible exception in the complex plane.
\end{lemma}

The next result, originally established by Hurwitz, is a foundational theorem that
plays a pivotal role in our analysis. It concerns the behavior of
zeros in sequences of analytic functions.

\begin{lemma}\label{R: Hurwitz}
If the sequence of non-zero analytic functions $f_n(z)$, defined on a domain $D$,
converges uniformly to a function $f$ on each compact subsets of $D$,
then $f(z)$ is either identically zero or never equal to zero in $D$.
\end{lemma}

Next, we recall the remarkable rescaling result due to
Lohwater and Pommerenke~\cite{Lohwater1973}, a result of central
importance in the study of normal meromorphic functions.

\begin{lemma}\cite{Lohwater1973}\label{R: Lohwater}
A function $f$, meromorphic in the unit disc $\mathbb{D}$, is a normal function
if and only if there do not exist sequences $\{z_n\}$ and $\{\rho_n\}$, with
$z_n \in \mathbb{D}$ and $\rho_n>0,\ \rho_n \to 0$, such that $g_n(z)=f(z_n+\rho_nz)$
converges uniformly on each compact subset of the complex plane to a non-constant meromorphic function $g(z)$.
\end{lemma}

With these preliminaries, we now proceed to the proof of our main theorems.

\section{Proof of Main results}

\begin{proof}[\bf{\textit{Proof of Theorem \ref{main}}}]
Assume that $f$, $\psi_1$, $\psi_2$, and $\psi_3$ satisfy the hypothesis
of the theorem and $f$ is not a normal function. Then, by Lemma~\ref{R: Lohwater},
there exist sequences $\{z_n\}$ and $\{\rho_n\}$, with $z_n \in \mathbb{D} $
and $\rho_n>0$, $\rho_n \rightarrow 0$ such that the sequence $\{g_n(z)=f(z_n+\rho_nz)\}$
converges uniformly on each compact subset of the complex plane
to a function $g(z)$, where $g(z)$ is a non-constant meromorphic function.
By taking a subsequence, if necessary, we may assume that
$z_n \to z_0 \in \overline{\mathbb{D}}$, the closure of $\mathbb{D}$.
Therefore, $z_n+\rho_nz \to z_0$ for each complex number $z$.\smallskip

If $z_0 \in \mathbb{D}, $ then $g_n(z)=f(z_n+\rho_nz)\to f(z_0) $,
which would mean $g(z)\equiv f(z_0)$, violating the assumption that $g$
is a non-constant function. Thus, we must have
that $|z_0|=1$, i.e. $z \in \partial \mathbb{D}$.\smallskip

Since $\psi_1$, $\psi_2$, and $\psi_3$ share the function $f$, thus we can consider the set,
\begin{equation}\label{1}
A=\{z\in \mathbb{D} :\psi_1(z)=f(z)\}=\{z\in \mathbb{D} :\psi_2(z)=f(z)\}=\{z\in \mathbb{D} : \psi_3(z)=f(z)\}.
\end{equation}
\textbf{Claim :} There are only finite terms of the sequence $\{z_n+\rho_nz\}$ will be in $A$.
Suppose, if our Claim holds true, then there exists $N \in \mathbb{N}$,
such that $z_n+\rho_n z \in A^c$ for all  $z \in \mathbb{C}$,
for all $n>N$. \smallskip

\noindent{Now, we consider two cases to prove our claim:}\smallskip

\noindent\textit{Case 1:} When $z_0 \not \in \overline{A}$.\smallskip

Then, we must get $n_0 \in \mathbb{N}$, such that $z_n+\rho_n z \in A^c $,
for all $n>n_0$, otherwise $z_0 \in \overline{A}$.\smallskip

\noindent\textit{Case 2:}  When $z_0\in \overline A$,
implies $z_0\in \overline{A}\cap\partial \mathbb{D}$
(since $z_0 \in \partial \mathbb{D} $).\smallskip

Suppose, on the contrary, that infinitely many terms of the sequence
$\{z_n+\rho_nz\}$ are in $A$, then we can find a subsequence
say $\{z_{n_k}+\rho_{n_k}z\}\in A$.\smallskip

Consider if possible, $f$ has pole for infinitely many
$ \zeta \in \{z_{n_k}+\rho_{n_k}z\}$ i.e. $f(\zeta)=\infty$
implies $\psi_i(\zeta)=\infty$, for all $i=1,2,3$.\smallskip

Since there are infinitely many such terms of sequence
$\{z_{n_k}+\rho_{n_k}z\}$ in $A$, so we can find a subsequence
of $\{z_{n_k}+\rho_{n_k}z \}$\,---\,relabeling if necessary, say
$\{z_{n_k}+\rho_{n_k}z \}$ such that,
$\psi_i(z_{n_k}+\rho_{n_k}z)=\infty$,
for all $i=1,2,3$, this implies
\begin{equation}\label{E: contra}
\frac{1}{\psi_i(z_{n_k}+\rho_{n_k}z)}=0.
\end{equation}

Since $\psi_i$ is continuous on $\partial{\mathbb{D}}$ for all $i=1,2,3$
and \[ \frac{1}{\psi_i(z_{n_k}+\rho_{n_k}z)} \rightarrow \frac{1}{\psi_i(z_0)} \Rightarrow
\frac{1}{\psi_i(z_0)}=0,\]
that is $\psi_i(z_0)=\infty$, for all $i=1,2,3$, contradiction to
the hypothesis that $\psi_i(z)\neq\psi_j(z)$ on $\partial\mathbb{D}$, $1\leq i<j\leq 3$.
Hence, there are only finitely many poles of $f$ in $\{z_{n_k}+\rho_{n_k}z\}$. Therefore,
we can find a subsequence of $\{z_{n_k}+\rho_{n_k}z\}$\,---\,relabelling if necessary,
say  $A \ni z_{n_k}+\rho_{n_k}z \rightarrow z_0 $
(ignoring the poles of $f$).\smallskip

Now fix $i, 1\leq i\leq 3$, and consider the function
$h_{n_k}(z)=f(z_{n_k}+\rho_{n_k}z)-\psi_i(z_{n_k}+\rho_{n_k}z)$,
thus $h_{n_k}\rightarrow g(z)-\psi_i(z_0)$ uniformly on each compact
subset of the plane, but $h_{n_k}(z_{n_k}+\rho_{n_k}z)=0$ for all $n_k
\in \mathbb{N}$ as $z_{n_k}+\rho_{n_k}z \in A$,
this implies $g(z)=\psi_i(z_0)$ for all $i=1,2,3$,
implies $\psi_1(z_0)=\psi_2(z_0)=\psi_3(z_0)$ which is a contradiction
to the hypothesis that $\psi_i(z)\neq\psi_j(z)$ on $\partial\mathbb{D}$, $1\leq i<j\leq 3$.
So, our assumption that infinitely many terms of the sequence
$\{z_n+\rho_nz\}$ are in $A$, is wrong.\smallskip

Hence our claim holds true, that there must exists
$n_1 \in \mathbb{N}$ such that $z_n+\rho_n z \in A^c$, for all $n>n_1, \ \ z \in \mathbb{C}$.\smallskip

Now in both the cases we will get a tail of sequence $z_n+\rho_n z $
which lies in $A^c$, by leaving first  $N=\max\{n_0,n_1\}$ terms.
Hence, our Claim holds true.\smallskip

Renumbering the sequence if necessary, without loss of generality
we may assume, $A^c \ni z_n+\rho_nz \to z_0 $ for all $z \in \mathbb{C} $.\smallskip

Fix $i, 1\leq i \leq 3$, and assume that $\psi_i(z_0)\neq \infty$,
$h_n(z)=f(z_n+\rho_n z)-\psi_i(z_n+\rho_nz)$, thus $h_n(z) \to g(z)-\psi_i(z)$
uniformly on each compact subset of the plane, since $f(z_n+\rho_nz)-\psi_i(z_n+\rho_nz)$
is never zero; (as $z_n+\rho_nz \in A^c$). It follows from Lemma \ref{R: Hurwitz}
due to Hurwitz, that either $g(z)-\psi_i(z_0)\equiv0$ or $g(z)-\psi_i(z_0)$ is
never zero, but $g(z)-\psi_i(z_0)\equiv0$ means that $g(z)$ is a constant
function, contradicting the assumption that it is not, therefore
it follows $g(z)$ never assumes the value $\psi_i(z_0)$.\smallskip

If $\psi_i(z_0)=\infty$, then we can take,
\[h_n^*(z)=\frac{1}{f(z_n+\rho_nz)}-\frac{1}{\psi_i(z_n+\rho_nz)}.
\]
And use the same argument (replacing $h^*$ in place of $h$ ) to conclude
that $1/g(z)$ does not assume the value 0. Which implies that $g(z)$
does not assume the value $\psi_i(z_0)=\infty$.\smallskip

Applying the same argument to each $i, 1\leq i \leq 3$, and using
the fact that $\psi_i(z)\neq\psi_j(z)$ $(1\leq i<j\leq 3)$  on $\partial\mathbb{D}$,
we get that the non-constant meromorphic function $g(z)$
avoids the three distinct values $\psi_1(z_0)$, $\psi_2(z_0)$, and $\psi_3(z_0)$.
But this contradicts Picard's Theorem. Thus, our assumption that $f$ is
not normal is unsound, hence $f$ is normal.\end{proof}\medskip

\begin{proof}[\bf{\textit{Proof of Theorem~\ref{submain}}}]
The outline of the proof remains the same as that of Theorem \ref{main}.
Again, we can assume that $f$ is not normal, then by Lemma~\ref{R: Lohwater}, there will be sequences
$\{z_n\}$ and $\{\rho_n\}$, with $z_n \in \mathbb{D} $ and
$\rho_n>0$, $\rho_n \rightarrow 0$ such that the sequence
$\{g_n(z)=f(z_n+\rho_nz)\}$ converges uniformly on each compact
subset of the complex plane to a function $g(z)$, where $g$ is a
non-constant meromorphic function. By taking a subsequence, if
necessary that $z_n \to z_0 \in \overline{\mathbb{D}}$, the
closure of $\mathbb{D}$. As done in Theorem \ref{main}, it is clear that $|z_0|=1$.
Since $\psi_1$ and $\psi_2$ share the function $f$,
thus we can consider the set,
\begin{equation}\label{2}
A=\{z\in \mathbb{D} :\psi_1(z)=f(z)\}=\{z\in \mathbb{D} :\psi_2(z)=f(z)\}
\end{equation}
We may proceed in the similar fashion as we have done in Theorem~\ref{main} and
show that there are only finite terms of the sequence $\{z_n+\rho_nz\}$ will be in $A$.\smallskip

we again consider two cases:\smallskip

\noindent\textit{Case 1:} If $z_0 \not \in \overline{A}$ \smallskip

Then, we must get $n_2 \in \mathbb{N}$ such that $z_n+\rho_n z \in A^c $ for all $n>n_2$.\smallskip

\noindent\textit{Case 2:} If $z_0\in \overline A$, implies $z_0\in \overline{A}\cap\partial \mathbb{D}$.\smallskip

Again suppose, if possible, that there are infinitely many terms of the sequence
$\{z_n+\rho_nz\}$ are in $A$, then we can find a subsequence
say $\{z_{n_k}+\rho_{n_k}z\}\in A$.\smallskip

Consider if possible, $f$ has pole for infinitely many
$ \zeta \in \{z_{n_k}+\rho_{n_k}z\}$ i.e. $f(\zeta)=\infty$,
implies $\psi_i(\zeta)=\infty$, for all $i=1,2$. Likewise, in the
\textit{Case 2} of the previous Theorem~\ref{main}, this time we will
get $\psi_1(z_0)=\psi_2(z_0)=\infty$. Which violates the assumption that
$\psi_1(z)\neq \psi_2(z)$ on $\partial\mathbb{D}$.\smallskip

Since, there are only finitely many points in the sequence
$\{z_{n_k}+\rho_{n_k}z\}$ where $f$ has poles, so we can find a
subsequence of $\{z_{n_k}+\rho_{n_k}z\}$\,---\, relabeling if necessary say $A \ni z_{n_k}+\rho_{n_k}z
\rightarrow z_0 $, (ignoring the poles of $f$).\smallskip

Now fix $i, 1\leq i\leq 2$, and consider the function
$h_{n_k}(z)=f(z_{n_k}+\rho_{n_k}z)- \psi_i(z_{n_k}+\rho_{n_k}z)$,
thus $h_{n_k}\rightarrow g(z)-\psi_i(z_0)$ uniformly on each compact
subset of the complex plane, but $h_{n_k}(z_{n_k}+\rho_{n_k}z)=0$ for all $n_k \in \mathbb{N}$
as $z_{n_k}+\rho_{n_k}z \in A$, this implies $g(z)=\psi_i(z_0)$  for all $i=1,2$,
implies $\psi_1(z_0)=\psi_2(z_0)$ which is a contradiction to the
hypothesis that $\psi_i(z)\neq\psi_j(z)$ on $\partial\mathbb{D}$ and $1\leq i<j\leq 2$.
So, our assumption that infinitely many terms of the sequence
$\{z_n+\rho_nz\}$ are in $A$, is wrong.\smallskip

Hence our claim holds true, that
there must exists $n_3 \in \mathbb{N}$ such that $z_n+\rho_n z \in A^c$
for all $n>n_3, z \in \mathbb{C}$.\medskip

Now from both the cases one can get a tail of sequence $z_n+\rho_n z $
which lies in $A^c$, by leaving first $N_{0}=\max\{n_2,n_3\}$ terms. \smallskip

Now proceeding with the similar manner as done in previous Theorem \ref{main}
and using the fact that $f$ avoids $\psi_3$ on $\mathbb{D}$,
we will get $g(z)$ avoids the three distinct values
$\psi_1(z_0)$, $\psi_2(z_0)$, and $\psi_3(z_0)$. A contradiction to the
fact that $g(z)$ is non-constant meromorphic function.
Hence, $f$ is normal.\end{proof}

\end{document}